\documentclass[a4paper,oneside,10pt]{article}%
\usepackage{amsmath}
\usepackage{amsfonts}
\usepackage{amssymb}
\usepackage{graphicx}
\usepackage[square,numbers,sort&compress]{natbib}%
\usepackage[toc,page,title,titletoc,header]{appendix}
\DeclareMathOperator*{\esssup}{ess\,sup}
\DeclareMathOperator*{\essinf}{ess\,inf}
\providecommand{\U}[1]{\protect\rule{.1in}{.1in}}

\newtheorem{theorem}{Theorem}[section]

\newtheorem{corollary}[theorem]{Corollary}

\newtheorem{definition}[theorem]{Definition}
\newtheorem{assumption}[theorem]{Assumption}
\newtheorem{example}[theorem]{Example}

\newtheorem{lemma}[theorem]{Lemma}

\newtheorem{proposition}[theorem]{Proposition}

\numberwithin{equation}{section}

\begin{document}

\title{The minimum mean square estimator of integrable variables under sublinear operators}
\author{Shaolin Ji\thanks{Zhongtai Institute of Finance, Shandong University, Jinan,
Shandong 250100, PR China. jsl@sdu.edu.cn. This research is supported by National Natural Science Foundation of China (No. 11571203), the Programme of Introducing Talents of Discipline to Universities of China (No. B12023).} \quad Chuiliu Kong\thanks{Corresponding author. Zhongtai Institute of Finance, Shandong University, Jinan,
Shandong 250100, PR China. kclsdmath@mail.sdu.edu.cn.}\quad Chuanfeng Sun\thanks{School of Mathematical Sciences, University of Jinan, Jinan, Shandong 250022, P.R. China. sms\_suncf@ujn.edu.cn. This research is partially supported by the National Natural Science Foundation of China (No. 11701214), the Natural Science Foundation of Shandong Province (No. ZR2017BA032).}}
\date{}
\maketitle
\textbf{Abstract}. In this paper, we study the minimum mean square estimator for non-bounded random variables under sublinear operators. The existence and uniqueness of the minimum mean square estimator are obtained. Several properties of the minimum mean square estimator for non-bounded random variables are proved  under some mild assumptions.

\section{Introduction}
In recent decades, nonlinear risk measures and nonlinear expectations have been proposed and developed rapidly. For example, Artzner et al.\cite{ADEHK} introduced coherent risk measure theory; Peng studied \emph{g}-expectation in \cite{P} and Related conditional nonlinear expectations have also been proposed.

It is well-known that, for the classical linear expectation case, the conditional expectation coincide with the minimum mean square estimator. From another view point, the minimum mean square estimator can also be used as an alternative definition of the conditional expectation. However, for the nonlinear expectation cases, we do not know the relation between the conditional nonlinear expectations and the minimum mean square estimator. Recently, Ji and Sun \cite{JS} introduced a new conditional nonlinear expectation for bounded random variables which is based on the minimum mean square estimator for sublinear operators. In their paper, they proved the existence and uniqueness of the minimum mean square estimator and give the basic properties of the minimum mean square estimator. The relationship between the minimum mean square estimator and the conditional coherent risk measure and conditional \emph{g}-expectation was explored.

However, the boundedness assumption for random variables in  \cite{JS} has great limitations. In this paper, our goal is to delete the boundedness assumption in  \cite{JS} and generalize the corresponding results to the case in which the random variables fall in the space $L_\mathcal{F}^{2+\epsilon}(\Omega, P_0)$ where  $\epsilon$ is a constant such that $\epsilon\in(0,1)$. To solve the minimum mean square estimate problem, we formulate it as a minimax problem due to that the sublinear operator can be represented as a supremum of a family of linear expectations. In more details, for the existence result we prove Proposition \ref{random variable limited} and Proposition \ref{P exist} which are necessary to construct a sequence of the optimal estimators in bounded integrable spaces $L_\mathcal{F}^{2+\epsilon,M}(\Omega, P_0)$ where $M<\infty$ is a constant. Based on the existence result, we obtain the form of the optimal estimator by the minimax theorem and use a construction method to deduce the uniqueness result. Comparing with some fundamental properties of the classical linear expectation, we prove that these properties for the minimum mean square estimators are also reasonable. At last, we illustrate the differences among the minimum square estimator, the conditional coherent risk measure and the conditional $\emph{g}$-expectation by three examples.

This paper is organized as follows. In section 2, we give some basic definitions and results and formulate our problem. In section 3, under some mild assumptions, the existence and uniqueness of the optimal estimator are established. In the last section, we prove the basic properties of the minimum mean square estimator and also explore the relationship between the minimum mean square estimator and the conditional coherent risk measure and conditional \emph{g}-expectation.

\section{Preliminary}
For a given complete probability space $(\Omega, \mathcal{F}, P_0)$, we denote the class of all $\mathcal{F}$-measurable $(2+\epsilon)$ integrable random variables by $L_\mathcal{F}^{2+\epsilon}(\Omega, P_0)$. Sometimes $L_\mathcal{F}^{2+\epsilon}(P_0)$ for short.
\begin{definition}
A sublinear operator is an operator $\rho:L_\mathcal{F}^{2+\epsilon}(\Omega, P_0)\rightarrow \mathbb{R}$ satisfying\\
(i)Monotonicity: $\rho(\xi_1)\geq \rho(\xi_2)$ \ if $\xi_1\geq\xi_2$;\\
(ii)Constant preserving: $\rho(c)=c$ \ for $c\in R$;\\
(iii)Sub-additivity: For each $\xi_1, \xi_2\in L_\mathcal{F}^{2+\epsilon}(\Omega, P_0)$, $\rho(\xi_1+\xi_2)\leq \rho(\xi_1)+\rho(\xi_2)$;\\
(iv)Positive homogeneity: $\rho(\lambda\xi)=\lambda\rho(\xi)$ for every constant $\lambda\geq 0$.
\end{definition}

\begin{theorem}\label{expression}
If $\rho$ is a sublinear operator and $\mathcal{P}$ is the family of all linear operators dominated by $\rho$, then
$$\rho(\xi)=\max_{P\in\mathcal{P}}E_P[\xi],\quad \forall\xi\in L_{\mathcal{F}}^{2+\epsilon}(\Omega,P_0).$$
\end{theorem}
\textbf{Proof.} By Corollary 2.4 of Chapter I in \cite{S}, for any $\xi\in L_{\mathcal{F}}^{2+\epsilon}(\Omega,P_0)$, there exists a linear operator $L$ such that $L<\rho$ and $L(\xi)=\rho(\xi)$. If we take all linear expectations dominated by $\rho$, then
$$\rho(\xi)=\max_{P\in\mathcal{P}}E_P[\xi],\quad \forall\xi\in L_{\mathcal{F}}^{2+\epsilon}(\Omega,P_0).$$
\rightline{$\square$}

Note that $L_\mathcal{F}^{2+\epsilon}(\Omega, P_0)$ is a Reflexive space. Denote the dual space of $L_\mathcal{F}^{2+\epsilon}(\Omega, P_0)$ by $\big(L_\mathcal{F}^{2+\epsilon}(\Omega, P_0)\big)^{*}$. By Theorem \ref{expression}, $\rho$ can be represented by the family of linear operators dominated by $\rho$. We also denote by $\mathcal{P}$ all linear operators dominated by $\rho$. The set $\mathcal{P}$ is called the representation set of $\rho$.

We need the following two assumptions. Unless indicated, this two assumptions are required throughout this paper.

\begin{assumption}\label{assumption proper}
The sublinear operator $\rho$ is proper, that is, all the elements in $\mathcal{P}$ are equivalent to $P_{0}$. Recall that two probability measures
$P$ and $P_{0}$ are said to be equivalent if for $A\in\mathcal{F}$, $P(A)=0$ if and only if $P_{0}(A)=0$.
\end{assumption}

\begin{assumption}\label{assumption compactness}
$\mathcal{D}:=\{\frac{dP}{dP_0},P\in\mathcal{P}\}$ is normed uniformly bounded in $L^{1+\frac{2}{\epsilon}}_{\mathcal{F}}(P_0)$ and $\sigma(L^{1+\frac{2}{\epsilon}}(P_0),L^{1+\frac{\epsilon}{2}}(P_0))$-compact.
\end{assumption}


In the sequel, for convenience, we will use $f^P$ to denote the Radon-Nikodym derivative $\frac{dP}{dP_0}$.
\begin{definition}[Stability]
We say that the set $\mathcal{P}$ is stable, if for each element $P\in\mathcal{P}$ with associated $f^{P}_{\mathcal{C}}, g_{\mathcal{C}}=\frac{f^{P}}{f^{P}_{\mathcal{C}}}$ still lies in $\mathcal{D}$, where $f^{P}_{\mathcal{C}}:=E_{P_0}[\frac{dP}{dP_0}|\mathcal{C}]$ and $\mathcal{C}$ is a sub-$\sigma$-algebra of $\mathcal{F}$ .
\end{definition}

We call the sublinear operator $\rho$ is stable, if its representation $\mathcal{P}$ is stable.
\begin{proposition}\label{stable}
If a sublinear operator $\rho$ is stable and proper, then for any $P\in\mathcal{P}$ and $\xi$ which is a integrable random variable, there exists a $\bar{P}\in\mathcal{P}$ such that $E_{\bar{P}}[\xi]=E_{P_0}[E_P[\xi|\mathcal{C}]].$
\end{proposition}
$\textbf{Proof}$. Since
$$E_{P_0}[E_P[\xi|\mathcal{C}]]=E_{P_0}\big[\frac{E_{P_0}[\xi f^P|\mathcal{C}]}{E_{P_0}[ f^P|\mathcal{C}]}\big]=E_{P_0}\big[E_{P_0}[\xi\frac{f^P}{f^P_{\mathcal{C}}}|\mathcal{C}]\big]=E_{P_0}[\xi\frac{f^P}{f^P_{\mathcal{C}}}].$$
Because of $\rho$ is stable, there exists $\bar{P}\in \mathcal{P}$ such that $\frac{d\bar{P}}{dP_0}=\frac{f^P}{f^P_{\mathcal{C}}}$. This implies  $E_{\bar{P}}[\xi]=E_{P_0}[E_P[\xi|\mathcal{C}]]$.\\
\rightline{$\square$}

Let $\mathcal{C}$ be a sub-$\sigma$-algebra of $\mathcal{F}$. For a given $\xi\in L^{4+2\epsilon}_{\mathcal{F}}(P_0)$, our problem is to find its minimum square estimator for the sublinear operator $\rho$ when ``the only information $\mathcal{C}$" is known for us, that is, to solve the following optimization problem.\\
$\textbf{Problem}$ \quad Find a $\hat{\eta}\in L^{2+\epsilon}_{\mathcal{C}}(P_0)$ such that
\begin{equation}\label{problem}
\rho(\xi-\hat{\eta})^2=\inf_{\eta\in L^{2+\epsilon}_{\mathcal{C}}(P_0)}\rho(\xi-\eta)^2
\end{equation}

The optimal solution $\hat{\eta}$ of \eqref{problem} is called the minimum mean square estimator. It is also regarded as a minimax estimator in statistical decision theory.

\section{Existence and Uniqueness Results}
In this section, we will study the existence and uniqueness of the minimum mean square estimator.
\subsection{Existence Result}
\begin{lemma}\label{le1}
For $\xi\in L^{2+\epsilon}_{\mathcal{F}}(\Omega,P_0)$, we have $\sup\limits_{P\in\mathcal{P}}E_P[\xi^2]<\infty.$
\end{lemma}
$\textbf{Proof}$. Since $\{f^P:P\in\mathcal{P}\}\in \mathcal{D}$ is normed uniformly bounded in $L^{1+\frac{2}{\epsilon}}_{\mathcal{F}}(P_0)$, it results
$$\sup_{P\in\mathcal{P}}E_P[\xi^2]=\sup_{P\in\mathcal{P}}E_{P_0}[f^P\xi^2]\leq \sup_{P\in\mathcal{P}}\|f^P\|_{L^q_{\mathcal{F}}}\|\xi^2\|_{L^p_{\mathcal{F}}}<\infty$$
where $p=\frac{2+\epsilon}{2}\ \mbox{and}\ q=\frac{2+\epsilon}{\epsilon}.$\\
\rightline{$\square$}
\begin{proposition}\label{random variable limited}
If $\xi\in L^{4+2\epsilon}_{\mathcal{F}}(\Omega,P_0)$ and the sublinear operator $\rho$ is stable, then there exists a constant $M$ such that for any $P\in\mathcal{P}$
$$\inf_{\eta\in L^{2+\epsilon}_{\mathcal{C}}(P_0)}E_P[(\xi-\eta)^2]=\inf_{\eta\in L^{2+\epsilon,M}_{\mathcal{C}}(P_0)}E_P[(\xi-\eta)^2]$$
where $L^{2+\epsilon,M}_{\mathcal{C}}(P_0)$ denotes all the elements in $L^{2+\epsilon}_{\mathcal{C}}(P_0)$ normed bounded by constant $M$.
\end{proposition}

$\textbf{Proof}$. Denote $\mathbb{G}:=\{E_P[\xi|\mathcal{C}];P\in\mathcal{P}\}$. For any $P\in\mathcal{P}$, the following relations hold
\begin{equation*}
\begin{aligned}
E_{P_0}[(E_P[\xi|\mathcal{C}])^{2+\epsilon}]&= E_{P_0}\big[\big(E_P[\xi\big|\mathcal{C}]^2\big)^{\frac{2+\epsilon}{2}}\big]\leq E_{P_0}\big[\big(E_P[\xi^2\big|\mathcal{C}]\big)^{\frac{2+\epsilon}{2}}\big]\\
&\leq E_{P_0}\big[E_{P}[\xi^{2\cdot\frac{2+\epsilon}{2}}\big|\mathcal{C}]\big]=E_{P_0}\big[E_P[\xi^{2+\epsilon}\big|\mathcal{C}]\big]
\end{aligned}
\end{equation*}
where the second $'\leq'$ comes from Jensen's inequality and the function $(x)^{1+\frac{\epsilon}{2}}$ is convex about $x$ when $x\geq 0$. By Proposition \ref{stable}, there exists a $\bar{P}\in \mathcal{P}$ such that $E_{\bar{P}}[\xi^{2+\epsilon}]=E_{P_0}[E_P[\xi^{2+\epsilon}|\mathcal{C}]]$. By Lemma \ref{le1}, there exists a constant $M_1$ such that $\sup_{P\in\mathcal{P}}E_P[\xi^{2+\epsilon}]\leq M_1$. Then  $\mathbb{G}\subset L^{{2+\epsilon},M}_{\mathcal{C}}(P_0)$, where $M=M_1^{\frac{1}{{2+\epsilon}}}$. Since
$$\mathbb{G}\subset L^{2+\epsilon,M}_{\mathcal{C}}(P_0)\subset L^{2+\epsilon}_{\mathcal{C}}(P_0)$$
and
$$E_P[(\xi-E_P[\xi|\mathcal{C}])^2]\leq E_P[(\xi-\eta)^2],\quad \forall \eta\in L^{2+\epsilon}_{\mathcal{C}}(P_0),$$
it results
$$\sup_{P\in \mathcal{P}}\inf_{\eta\in L^{2+\epsilon}_{\mathcal{C}}(P_0)}\big[E_P[(\xi-\eta)^2]\big]\geq\sup_{P\in \mathcal{P}}\inf_{\eta'\in \mathbb{G}}\big[E_P[(\xi-\eta')^2]\big].$$
On the other hand, since $\mathbb{G}\subset L^{2+\epsilon}_{\mathcal{C}}(P_0)$, the inverse inequality is obviously true. Then the following equality holds
$$\sup_{P\in \mathcal{P}}\inf_{\eta\in L^{2+\epsilon}_{\mathcal{C}}(P_0)}\big[E_P[(\xi-\eta)^2]\big]=\sup_{P\in \mathcal{P}}\inf_{\eta'\in \mathbb{G}}\big[E_P[(\xi-\eta')^2]\big].$$
Hence, it follows that
$$\sup_{P\in \mathcal{P}}\inf_{\eta\in L^{2+\epsilon}_{\mathcal{C}}(P_0)}\big[E_P[(\xi-\eta)^2]\big]=\sup_{P\in \mathcal{P}}\inf_{\eta'\in L^{2+\epsilon,M}_{\mathcal{C}}(P_0)}\big[E_P[(\xi-\eta' )^2]\big].$$\\
\rightline{$\square$}
\begin{proposition}\label{P exist}
For a given $\xi\in L^{4+2\epsilon}_{\mathcal{F}}(\Omega,P_0)$, the following equality holds
$$\sup_{P\in \mathcal{P}}\inf_{\eta\in L^{2+\epsilon,M}_{\mathcal{C}}(P_0)}\big[E_P[(\xi-\eta)^2]\big]=\max_{P\in \mathcal{P}}\inf_{\eta\in L^{2+\epsilon,M}_{\mathcal{C}}(P_0)}\big[E_P[(\xi-\eta)^2]\big].$$
\end{proposition}
$\textbf{Proof}$. Let
$$\beta:=\sup_{P\in \mathcal{P}}\inf_{\eta\in L^{2+\epsilon,M}_{\mathcal{C}}(P_0)}\big[E_P[(\xi-\eta)^2]\big]=\sup_{f^{P}\in \mathcal{D}}\inf_{\eta\in L^{2+\epsilon,M}_{\mathcal{C}}(P_0)}\big[E_{P_0}[f^{P}(\xi-\eta)^2]\big].$$
Take a sequence $\{f^{P_n};P_n\in \mathcal{P}\}_{n\geq 1}$ such that
$$\inf_{\eta\in L^{2+\epsilon,M}_{\mathcal{C}}(P_0)}\big[E_{P_0}[f^{P_n}(\xi-\eta)^2]\big]\geq \beta-\frac{1}{2^n}.$$
Since the set $\mathcal{D}$ is a weakly compact set, we can take a subsequence$\{f^{P_{n_i}}\}_{i\geq 1}$ of $\{f^{P_n};P_n\in \mathcal{P}\}_{n\geq 1}$ which weakly converges to some $f^{\hat{P}}\in L^{1+\frac{2}{\epsilon}}(P_0)$.  Therefore, thanks to a separation Hahn-Banach standard result, there exists a sequence $\{f^{\tilde{P}}_i\in conv(f^{P_{n_i}},f^{P_{n_{i+1}}},...)\}_{i\geq 1}$ such that $f^{\tilde{P}}_i$ converges to $f^{\hat{P}}$ in $L^{1+\frac{2}{\epsilon}}(P_0)$-norm. This shows that $\hat{P}\in \mathcal{P}$.

On the other hand, for any $\eta\in L^{2+\epsilon,M}_{\mathcal{C}}(P_0)$ and $i\in \mathbb{N}$, the following inequality holds
$$E_{P_0}[f^{\tilde{P}}_i(\xi-\eta)^2]\geq \inf_{\tilde{\eta}\in L^{2+\epsilon,M}_{\mathcal{C}}(P_0)}E_{P_0}[f^{\tilde{P}}_i(\xi-\tilde{\eta})^2].$$
Then for any $\eta\in L^{2+\epsilon,M}_{\mathcal{C}}(P_0)$, it follows that
$$\lim_{i\rightarrow\infty}E_{P_0}[f^{\tilde{P}}_i(\xi-\eta)^2]\geq \limsup_{i\rightarrow\infty}\inf_{\tilde{\eta}\in L^{2+\epsilon,M}_{\mathcal{C}}(P_0)}E_{P_0}[f^{\tilde{P}}_i(\xi-\eta)^2].$$
Thus
\begin{equation}\label{p3.1}
\inf_{\eta\in L^{2+\epsilon,M}_{\mathcal{C}}(P_0)}\lim_{i\rightarrow\infty}E_{P_0}[f^{\tilde{P}}_i(\xi-\eta)^2]\geq \limsup_{i\rightarrow\infty}\inf_{\tilde{\eta}\in L^{2+\epsilon,M}_{\mathcal{C}}(P_0)}E_{P_0}[f^{\tilde{P}}_i(\xi-\eta)^2].
\end{equation}
Since $||(\xi-\eta)^2||_{L^{1+\frac{\epsilon}{2}}(P_0)}< \infty$, it results that
$$\lim_{i\rightarrow\infty}E_{P_0}|f^{\tilde{P}}_i(\xi-\eta)^2-f^{\hat{P}}(\xi-\eta)^2|\leq \lim_{i\rightarrow\infty}||(f^{\tilde{P}}_i-f^{\hat{P}})||_{L^{1+\frac{2}{\epsilon}}(P_0)}||(\xi-\eta)^2||_{L^{1+\frac{\epsilon}{2}}(P_0)}= 0.$$
Then
$$E_{P_0}[f^{\hat{P}}(\xi-\eta)^2]=\lim_{i\rightarrow\infty}E_{P_0}[f^{\tilde{P}}_{i}(\xi-\eta)^2].$$
It results
\begin{equation}\label{p3.2}
\inf_{\eta\in L^{2+\epsilon,M}_{\mathcal{C}}(P_0)}E_{P_0}[f^{\hat{P}}(\xi-\eta)^2]=\inf_{\eta\in L^{2+\epsilon,M}_{\mathcal{C}}(P_0)}\lim_{i\rightarrow\infty}E_{P_0}[f^{\tilde{P}}_i(\xi-\eta)^2].
\end{equation}
By \eqref{p3.1} and \eqref{p3.2}, the following relations hold
$$\inf_{\eta\in L^{2+\epsilon,M}_{\mathcal{C}}(P_0)}E_{P_0}[f^{\hat{P}}(\xi-\eta)^2]\geq \limsup_{i\rightarrow\infty}\inf_{\tilde{\eta}\in L^{2+\epsilon,M}_{\mathcal{C}}(P_0)}E_{P_0}[f^{\tilde{P}}_i(\xi-\eta)^2]\geq \beta.$$
Since $\hat{P}\in\mathcal{P}$, we get
$$\inf_{\eta\in L^{2+\epsilon,M}_{\mathcal{C}}(P_0)}E_{P_0}[f^{\hat{P}}(\xi-\eta)^2]=\sup_{P\in \mathcal{P}}\inf_{\eta\in L^{2+\epsilon,M}_{\mathcal{C}}(P_0)}\big[E_P[(\xi-\eta)^2]\big].$$\\
\rightline{$\square$}

\begin{corollary}\label{P exist2}
If the sublinear operator $\rho$ is stable, then for a given $\xi\in L^{4+2\epsilon}_{\mathcal{F}}(\Omega,P_0)$, the following equality holds
$$\sup_{P\in \mathcal{P}}\inf_{\eta\in L^{2+\epsilon}_{\mathcal{C}}(P_0)}E_P[(\xi-\eta)^2]=\max_{P\in \mathcal{P}}\inf_{\eta\in L^{2+\epsilon}_{\mathcal{C}}(P_0)}E_P[(\xi-\eta)^2].$$
\end{corollary}
$\textbf{Proof}$. Choose $\hat{P}$ as in Proposition \ref{P exist}. By Propositions \ref{random variable limited} and \ref{P exist}, the following relations hold
\begin{equation*}
\begin{aligned}
&\sup_{P\in \mathcal{P}}\inf_{\eta\in L^{2+\epsilon}_{\mathcal{C}}(P_0)}E_P[(\xi-\eta)^2]= \sup_{P\in \mathcal{P}}\inf_{\eta\in L^{2+\epsilon,M}_{\mathcal{C}}(P_0)}E_P[(\xi-\eta)^2]=\inf_{\eta\in L^{2+\epsilon,M}_{\mathcal{C}}(P_0)}E_{\hat{P}}[(\xi-\eta)^2]\\
&=\inf_{\eta\in L^{2+\epsilon}_{\mathcal{C}}(P_0)}E_{\hat{P}}[(\xi-\eta)^2].
\end{aligned}
\end{equation*}
Since $\hat{P}\in\mathcal{P}$, one obtains
$$\sup_{P\in \mathcal{P}}\inf_{\eta\in L^{2+\epsilon}_{\mathcal{C}}(P_0)}E_P[(\xi-\eta)^2]=\max_{P\in \mathcal{P}}\inf_{\eta\in L^{2+\epsilon}_{\mathcal{C}}(P_0)}E_P[(\xi-\eta)^2].$$
\rightline{$\square$}

\begin{theorem}[Fan.K \cite{FK}(1953)]\label{Minimax Theorem}
Let $\mathcal{X}$ be a compact Hausdorff space and $\mathcal{Y}$ be an arbitrary set. Let $F$ be a real valued function defined on $\mathcal{X}\times \mathcal{Y}$ such that, for every $y\in\mathcal{Y}$, $F(x,y)$ is a $l.s.c$(lower-semicontinuous) on $\mathcal{X}$. If $F$ is convex on $\mathcal{X}$ and concave on $\mathcal{Y}$, then\\
$$\min_{x\in\mathcal{X}}\sup_{y\in\mathcal{Y}}F(x,y)=\sup_{y\in\mathcal{Y}}\min_{x\in\mathcal{X}}F(x,y).$$
\end{theorem}
$\textbf{Proof}$. Refer to Theorem 2 in \cite{FK}.\\
\rightline{$\square$}

\begin{theorem}[Existence Theorem]\label{Existence Theorem}
If $\xi\in L^{4+2\epsilon}_{\mathcal{F}}(\Omega,P_0)$ and the sublinear operator $\rho$ is stable, then there exists an optimal solution $\hat{\eta}\in L^{2+\epsilon}_{\mathcal{C}}(\Omega,P_0)$ for the $\textbf{Problem}$ \eqref{problem}.
\end{theorem}
$\textbf{Proof}$. Since $\xi\in L^{4+2\epsilon}_{\mathcal{F}}(P_0)$ and $\eta\in L^{2+\epsilon}_{\mathcal{F}}(P_0)$, then
$(\xi-\eta)^2\in L^{1+\frac{\epsilon}{2}}_{\mathcal{F}}(P_0)$. This means $E_{P_0}[f^P(\xi-\eta)^2]$ is a continuous function on topology space $(L^{1+\frac{2}{\epsilon}}(P_0),\sigma(L^{1+\frac{2}{\epsilon}}(P_0),L^{1+\frac{\epsilon}{2}}(P_0)))$. Because the set $\mathcal{D}$ is $\sigma(L^{1+\frac{2}{\epsilon}}(P_0),L^{1+\frac{\epsilon}{2}}(P_0))$-compact, then by Theorem \ref{Minimax Theorem} and Corollary \ref{P exist2}, the following equality holds
$$\max_{P\in \mathcal{P}}\inf_{\eta\in L^{2+\epsilon}_{\mathcal{C}}(P_0)}E_P[(\xi-\eta)^2]=\inf_{\eta\in L^{2+\epsilon}_{\mathcal{C}}(P_0)}\max_{P\in \mathcal{P}}E_P[(\xi-\eta)^2].$$
Moreover, with the help of Proposition \ref{random variable limited}, we derive
$$\max_{P\in \mathcal{P}}\inf_{\eta'\in L^{2+\epsilon,M}_{\mathcal{C}}(P_0)}E_P[(\xi-\eta' )^2]=\inf_{\eta'\in L^{2+\epsilon,M}_{\mathcal{C}}(P_0)}\max_{P\in \mathcal{P}}E_P[(\xi-\eta' )^2].$$
Therefore, it results
$$\inf_{\eta'\in L^{2+\epsilon,M}_{\mathcal{C}}(P_0)}\max_{P\in \mathcal{P}}E_P[(\xi-\eta' )^2]=\inf_{\eta\in L^{2+\epsilon}_{\mathcal{C}}(P_0)}\max_{P\in \mathcal{P}}E_P[(\xi-\eta)^2].$$
Hence, we can take a sequence $\{\eta_n;n\in \mathbb{N}\}\subset  L^{2+\epsilon,M}_{\mathcal{C}}(P_0)$ such that
$$\rho(\xi-\eta_n)^2<\alpha+\frac{1}{2^n}$$
where $\alpha:=\inf\limits_{\eta\in L^{2+\epsilon}_{\mathcal{C}}(P_0)}\rho(\xi-\eta)^2$. Since $L^{2+\epsilon,M}_{\mathcal{C}}(P_0)$ is a weakly compact set, we can take a subsequence $\{\eta_{n_i}\}_{i\in \mathbb{N}}$ of $\{\eta_n\}_{n\in \mathbb{N}}$ which weakly converges to some $\hat{\eta}\in L^{2+\epsilon,M}_{\mathcal{C}}(P_0)$. Using a separation of convex sets Hahn-Banach result, there exists a sequence $\{\tilde{\eta}_i\in conv(\eta_{n_i},\eta_{n_{i+1}},...)\}_{i\in \mathbb{N}}$ such that $\tilde{\eta}_i$ converges to $\hat{\eta}$ in $L^{2+\epsilon}_{\mathcal{C}}(P_0)$-norm. Since
\begin{equation}\label{Calculate}
\begin{aligned}
\rho(\xi-\hat{\eta})^2&=\rho(\xi-\tilde{\eta}_i+\tilde{\eta}_i-\hat{\eta})^2\\
&=\sup_{P\in\mathcal{P}}E_P[(\xi-\tilde{\eta}_i)^2+(\tilde{\eta}_i-\hat{\eta})^2+2(\xi-\tilde{\eta}_i)(\tilde{\eta}_i-\hat{\eta})]\\
&\leq \sup_{P\in\mathcal{P}}E_P[(\xi-\tilde{\eta}_i)^2]+\sup_{P\in\mathcal{P}}E_P[(\tilde{\eta}_i-\hat{\eta})^2+2(\xi-\tilde{\eta}_i)(\tilde{\eta}_i-\hat{\eta})]\\
&=\rho(\xi-\tilde{\eta}_i)^2+\sup_{P\in\mathcal{P}}E_P[-(\tilde{\eta}_i-\hat{\eta})^2+2(\xi-\hat{\eta})(\tilde{\eta}_i-\hat{\eta})]\\
&\leq \alpha+\frac{1}{2^{i-1}}+2\sup_{P\in\mathcal{P}}\|f^P\|_{L^{1+\frac{2}{\epsilon}}(P_0)}\|(\xi-\hat{\eta})\|_{L^{2+\epsilon}(P_0)}\|(\tilde{\eta}_i-\hat{\eta})\|_{L^{2+\epsilon}(P_0)}.\\
\end{aligned}
\end{equation}
Since \eqref{Calculate} holds for any $i\geq 1$, one obtains $\rho(\xi-\hat{\eta})^2=\alpha$.\\
\rightline{$\square$}

\subsection{Uniqueness Theorem}
In this sequel, we prove that the optimal solution of $\textbf{Problem}$ \eqref{problem} is unique.
\begin{theorem}[$\breve{\mathbf{Z}}$alinescu \cite{Z} (2002)]\label{saddle point}
Let A and B be two nonempty sets and $f$ from A$\times$ B to $\mathbb{R}\cup \{\infty\}$. Then $f$ has saddle points,i.e.there exists $(\bar{x},\bar{y})\in A\times B$, such that
$$\forall x\in A,\ \forall y\in B:\quad f(x,\bar{y})\leq f(\bar{x},\bar{y})\leq f(\bar{x},y)$$
if and only if
$$\inf_{y\in B}f(\bar{x},y)=\max_{x\in A}\inf_{y\in B}f(x,y)=\min_{y\in B}\sup_{x\in A}f(x,y)=\sup_{x\in A}f(x,\bar{y}).$$
\end{theorem}
$\textbf{Proof}$. Refer to Theorem 2.10.1 of Chapter 2 in \cite{Z}.\\
\rightline{$\square$}

\begin{theorem}\label{unique}
If the sublinear operator $\rho$ is stable, then the optimal solution of problem \eqref{problem} is unique.
\end{theorem}
$\textbf{Proof}$. From Theorem \ref{Existence Theorem}, the optimal solution exists. In the rest, we prove the optimal solution is unique. By Theorem \ref{Minimax Theorem} and Corollary \ref{P exist2}, the following equality holds
$$\max_{P\in \mathcal{P}}\inf_{\eta\in L^{2+\epsilon}_{\mathcal{C}}(P_0)}E_P[(\xi-\eta)^2]=\inf_{\eta\in L^{2+\epsilon}_{\mathcal{C}}(P_0)}\max_{P\in \mathcal{P}}E_P[(\xi-\eta)^2].$$
Since the optimal solution exists, it results
$$\max_{P\in \mathcal{P}}\min_{\eta\in L^{2+\epsilon}_{\mathcal{C}}(P_0)}E_P[(\xi-\eta)^2]=\min_{\eta\in L^{2+\epsilon}_{\mathcal{C}}(P_0)}\max_{P\in \mathcal{P}}E_P[(\xi-\eta)^2].$$
Denote the optimal solution by $\hat{\eta}$. By Corollary \ref{P exist2}, there exists $\hat{P}\in\mathcal{P}$ such that
$$\inf_{\eta\in L^{2+\epsilon}_{\mathcal{C}}(P_0)}E_{P_0}[f^{\hat{P}}(\xi-\eta)^2]=\max_{P\in \mathcal{P}}\inf_{\eta\in L^{2+\epsilon}_{\mathcal{C}}(P_0)}E_{P_0}[f^P(\xi-\eta)^2].$$
By Theorem \ref{saddle point}, the $(\hat{\eta},\hat{P})$ is the saddle point, i.e.
$$E_{P_0}[f^P(\xi-\hat{\eta})]^2\leq E_{P_0}[f^{\hat{P}}(\xi-\hat{\eta})]^2 \leq E_{P_0}[f^{\hat{P}}(\xi-\eta)]^2.$$
This shows that if $\hat{\eta}$ is the optimal solution, then there exists a $\hat{P}\in\mathcal{P}$ such that $\hat{\eta}=E_{\hat{P}}[\xi|\mathcal{C}]$.

Suppose that there exist two optimal solutions $\hat{\eta}_1$ and $\hat{\eta}_2$. Denote the accompanying probabilities by $\hat{P}_1$ and $\hat{P}_2$ respectively. Then we have $\hat{\eta}_1=E_{\hat{P}_1}[\xi|\mathcal{C}]$ and $\hat{\eta}_2=E_{\hat{P}_2}[\xi|\mathcal{C}]$. Set $P^{\lambda}=\lambda\hat{P}_1+(1-\lambda)\hat{P}_2$, $\lambda\in(0,1)$. Let $\lambda_{\hat{P}_1}=\lambda E_{P^{\lambda}}\big[\frac{d\hat{P}_1}{dP^{\lambda}}|\mathcal{C}\big]$ and $\lambda_{\hat{P}_2}=(1-\lambda )E_{P^{\lambda}}\big[\frac{d\hat{P}_2}{dP^{\lambda}}|\mathcal{C}\big]$ such that $\lambda_{\hat{P}_1}+\lambda_{\hat{P}_2}=1$. Then we have the following inequality (Details of the calculation can be found in Lemma \ref{calculate} in Appendix A):
\begin{equation}
\begin{aligned}
E_{P^{\lambda}}[(\xi-E_{P^{\lambda}}[\xi|\mathcal{C}])^{2}]=&E_{P^{\lambda}}[(\xi-\lambda_{\hat{P}_1}\hat{\eta}_1-\lambda_{\hat{P}_2}\hat{\eta}_2)^2]\\
=&E_{P^{\lambda}}\big[\big(\lambda_{\hat{P}_1}(\xi-\hat{\eta}_1)+\lambda_{\hat{P}_2}(\xi-\hat{\eta}_2)\big)^2\big]\\
=&E_{P^{\lambda}}\big[\lambda^2_{\hat{P}_1}(\xi-\hat{\eta}_1)^2+\lambda^2_{\hat{P}_2}(\xi-\hat{\eta}_2)^2+2\lambda_{\hat{P}_1}\lambda_{\hat{P}_2}(\xi-\hat{\eta}_1)(\xi-\hat{\eta}_1)\big]\\
=&E_{P^{\lambda}}\big[\lambda_{\hat{P}_1}(\xi-\hat{\eta}_1)^2+\lambda_{\hat{P}_2}(\xi-\hat{\eta}_2)^2-\lambda_{\hat{P}_1}\lambda_{\hat{P}_2}(\hat{\eta}_1-\hat{\eta}_2)^2\big]\\
=&\lambda E_{\hat{P}_1}\big[(\xi-\hat{\eta}_1)^2\big]+(1-\lambda) E_{\hat{P}_2}\big[(\xi-\hat{\eta}_2)^2\big]\\
&+\lambda E_{\hat{P}_1}\big[\lambda^2_{\hat{P}_2}(\hat{\eta}_1-\hat{\eta}_2)^2\big]+(1-\lambda)E_{\hat{P}_2}\big[\lambda^2_{\hat{P}_1}(\hat{\eta}_1-\hat{\eta}_2)^2\big]\\
\geq &\alpha
\end{aligned}
\end{equation}
where $\alpha:=\inf\limits_{\eta\in L^{2+\epsilon}_{\mathcal{C}}(P_0)}\rho(\xi-\eta)^2$.

Since $\rho$ is proper, then $E_{P^{\lambda}}[(\xi-E_{P^{\lambda}}[\xi|\mathcal{C}])^2]=\alpha$ if and only if $\hat{\eta}_1=\hat{\eta}_2$, $P_0$-a.s., i.e., $P_0(\{\omega: \hat{\eta}_{1}(\omega)=\hat{\eta}_{2}(\omega)\})=1$.

On the other hand, since $(\hat{\eta}_1,\hat{P}_1)$ is a saddle point, the following relations hold
$$E_{P^{\lambda}}\big[(\xi-E_{P^{\lambda}}[\xi|\mathcal{C}])^2\big]\leq E_{P^{\lambda}}\big[(\xi-\hat{\eta}_1)^2\big]\leq E_{\hat{P}_1}\big[(\xi-\hat{\eta}_1)^2\big]=\alpha.$$

It yields that $E_{P^{\lambda}}\big[(\xi-E_{P^{\lambda}}[\xi|\mathcal{C}])^2\big]=\alpha$. Thus, we deduce $\hat{\eta}_1=\hat{\eta}_2,$ $P_0$-a.s., i.e., $P_0(\{\omega: \hat{\eta}_{1}(\omega)=\hat{\eta}_{2}(\omega)\})=1.$\\
\rightline{$\square$}
\textbf {Remark.} We can also characterize the minimum mean square estimator like Ji and Sun in \cite{JS}, and give out the equivalent condition of optimal solution. So we omit this part in this paper.

\section{Properties of the Minimum Mean Square Estimator}
In this section, we will give some basic properties of the minimum mean square estimator. Then we explore the relationship between the minimum mean square estimator and the conditional coherent risk measure and conditional $\emph{g}$-expectation.

For a given $\xi\in L^{4+2\epsilon}_{\mathcal{F}}(P_0)$, we denote the minimum mean square estimator with respect to $\mathcal{C}$ by $\rho(\xi|\mathcal{C})$. Then $\rho(\xi|\mathcal{C})$ satisfies the following properties.
\begin{proposition}\label{proposition}
If the sublinear operator $\rho$ is stable and proper, then for any $\xi\in L^{4+2\epsilon}_{\mathcal{F}}(P_0)$, we obtain\\
i)If $C_1\leq \xi\leq C_2$ for two constants $C_1$ and $C_2$, then $C_1\leq \rho(\xi|\mathcal{C}) \leq C_2$.\\
ii)$\rho(\lambda\xi|\mathcal{C})=\lambda\rho(\xi|\mathcal{C})$ for any $\lambda\in \mathbb{R}$.\\
iii)For each $\eta_0\in L^{2+\epsilon}_{\mathcal{C}}(P_0)$, then $\rho(\xi+\eta_0|\mathcal{C})=\rho(\xi|\mathcal{C})+\eta_0$.\\
iv)If under each $P\in \mathcal{P}$, $\xi$ is independent of the sub $\sigma$-algebra $\mathcal{C}$, then $\rho(\xi|\mathcal{C})$ is a constant.
\end{proposition}
$\textbf{Proof}$. i) If $C_1\leq \xi\leq C_2$, then for any $P\in\mathcal{P}$, we have $C_1\leq E_P[\xi|\mathcal{C}]\leq C_2$. Since $\rho(\xi|\mathcal{C})\in \{E_P[\xi|\mathcal{C}]; P\in\mathcal{P}\}$, then $\rho(\xi|\mathcal{C})$ lies in $[C_1,C_2]$.\\
ii)If $\lambda=0$, the result is obvious. If $\lambda \neq 0$, it follows that
$$\lambda^2\rho(\xi-\frac{\rho(\lambda\xi|\mathcal{C})}{\lambda})^2=\rho(\lambda\xi-\rho(\lambda\xi|\mathcal{C}))^2=\inf_{\eta\in L^{2+\epsilon}_{\mathcal{C}}(P_0)}\rho(\lambda\xi-\eta)^2=\lambda^2\inf_{\eta\in L^{2+\epsilon}_{\mathcal{C}}(P_0)}\rho(\xi-\eta)^2.$$
It results that
$$\rho(\xi-\frac{\rho(\lambda\xi|\mathcal{C})}{\lambda})^2=\inf_{\eta\in L^{2+\epsilon}_{\mathcal{C}}(P_0)}\rho(\xi-\eta)^2.$$
Thus
$$\frac{\rho(\lambda\xi|\mathcal{C})}{\lambda}=\rho(\xi|\mathcal{C}).$$\\
iii) Note that
$$\rho(\xi+\eta_0-(\eta_0+\rho(\xi|\mathcal{C})))^2=\rho(\xi-\rho(\xi|\mathcal{C}))^2=\inf_{\eta\in L^{2+\epsilon}_{\mathcal{C}}(P_0)}\rho(\xi-\eta)^2=\inf_{\eta\in L^{2+\epsilon}_{\mathcal{C}}(P_0)}\rho(\xi+\eta_0-\eta)^2.$$
By the uniqueness of the minimum mean square estimator, the following equality holds
$$\rho(\xi+\eta_0|\mathcal{C})=\eta_0+\rho(\xi|\mathcal{C}).$$\\
iv) If under each $P\in \mathcal{P}$, $\xi$ is independent of the sub $\sigma$-algebra $\mathcal{C}$, then $E_P[\xi|\mathcal{C}]$ is a constant for each $P\in\mathcal{P}$. Since $\rho(\xi|\mathcal{C})\in \{E_P[\xi|\mathcal{C}]; P\in\mathcal{P}\}$, we know that $\rho(\xi|\mathcal{C})$ is a constant.\\
\rightline{$\square$}

The conditional coherent risk measure and some special conditional $\emph{g}$-expectations which were introduced by Artzner et al. \cite{ADEHK} and Peng \cite{P} respectively can be defined by $\esssup\limits_{P\in\mathcal{P}}E_P[\xi|\mathcal{C}]$. In the next three examples, we will show that the minimum mean square estimator is different from the conditional coherent risk measure and the conditional $\emph{g}$-expectation.
\begin{example}
Let $\Omega=\{\omega_1,\omega_2\}$, $\mathcal{F}=\{\phi,\{\omega_1\},\{\omega_2\},\Omega\}$ and $\mathcal{C}=\{\phi,\Omega\}$. Set $P_1=\frac{1}{3}I_{\omega_1}+\frac{2}{3}I_{\omega_2}$, $P_2=\frac{2}{3}I_{\omega_1}+\frac{1}{3}I_{\omega_2}$ and $\mathcal{P}=\{\lambda P_1+(1-\lambda)P_2; \lambda\in[0,1]\}$. For each $\xi\in L^{4+2\epsilon}_{\mathcal{F}}(P_0)$, define
$$\rho(\xi)=\sup_{P\in\mathcal{P}}E_P[\xi].$$
Set $\xi=2I_{\omega_1}+6I_{\omega_2}$. It is easy to see that
$$\sup_{P\in\mathcal{P}}E_P[\xi]=\frac{14}{3}\quad and \quad \rho(\xi|\mathcal{C})=E_{\hat{P}}[\xi|\mathcal{C}]=4$$
where $\hat{P}=\frac{1}{2}I_{\omega_1}+\frac{1}{2}I_{\omega_2}$.
\end{example}
\begin{example}
Let $\Omega=\{1,2,3,...\}$, $\mathcal{F}$ be the power set of $\Omega$ and $\mathcal{C}=\{\phi,\Omega\}$. Set
\begin{equation}
P_1=\left\{
\begin{aligned}
&\frac{1}{2},\quad\,\, \omega=1\\
&\frac{1}{2^{2}},\quad \omega=2\\
&\,\,\vdots \quad \quad \quad \vdots\\
&\frac{1}{2^{n}},\quad \omega=n\\
&\,\,\vdots \quad \quad \quad \vdots
\end{aligned}
\right.
\mbox{,}\quad
P_2=\left\{
\begin{aligned}
&\frac{2}{3},\quad\,\, \omega=1\\
&\frac{2}{3^{2}},\quad \omega=2\\
&\,\,\vdots \quad \quad \quad \vdots\\
&\frac{2}{3^{n}},\quad \omega=n\\
&\,\,\vdots \quad \quad \quad \vdots
\end{aligned}
\right.
\mbox{and}\quad
\xi=\left\{
\begin{aligned}
&1,\quad\,\, \omega=1\\
&\frac{1}{4},\quad \omega=2\\
&\,\,\vdots \quad \quad \quad \vdots\\
&\frac{2^{n}}{n^{4}},\quad \omega=n\\
&\,\,\vdots \quad \quad \quad \vdots
\end{aligned}
\right.
\end{equation}
and
\begin{equation}
\mathcal{P}=\{\lambda P_1+(1-\lambda)P_2; \lambda\in[0,1]\}=\left\{
\begin{aligned}
&\frac{\lambda}{2}+\frac{2(1-\lambda)}{3},\quad \omega=1,\\
&\frac{\lambda}{2^2}+\frac{2(1-\lambda)}{3^2},\quad \omega=2,\\
&\quad\vdots \quad \quad \quad \quad \quad \quad \quad\vdots\\
&\frac{\lambda}{2^n}+\frac{2(1-\lambda)}{3^n},\quad \omega=n,\\
&\quad\vdots \quad \quad \quad \quad \quad \quad \quad\vdots
\end{aligned}
\right.
\end{equation}
Define
$$\rho(\xi)=\sup_{P\in\mathcal{P}}E_P[\xi].$$
\end{example}
\begin{equation}
\begin{aligned}
E_P[\xi]&=[\frac{\lambda}{2}+\frac{2(1-\lambda)}{3}]+[\frac{\lambda}{2^2}+\frac{2(1-\lambda)}{3^2}]\cdot\frac{1}{2^2}+\cdots+[\frac{\lambda}{2^n}+\frac{2(1-\lambda)}{3^n}]\cdot\frac{2^n}{n^4}+\cdots\\
&=\frac{2}{3}+\sum_{n=2}^{\infty}\frac{2}{3^n}\cdot\frac{2^n}{n^4}+[-\frac{1}{6}+\sum_{n=2}^{\infty}(\frac{1}{n^4}-\frac{2^{n+1}}{3^n\cdot n^4})]\lambda.
\end{aligned}
\end{equation}
By $\sum_{n=1}^{\infty}\frac{1}{n^4}=\frac{\pi^4}{90}$, we have that $-\frac{1}{6}+\sum_{n=2}^{\infty}(\frac{1}{n^4}-\frac{2^{n+1}}{3^n\cdot n^4})<0$ which leads to
$$\sup_{P\in\mathcal{P}}E_P[\xi]=\frac{2}{3}+\sum_{n=2}^{\infty}\frac{2}{3^n}\cdot\frac{2^n}{n^4}.$$

Then, we calculate the optimal mean square estimator. For $p_n\geq0,\ n\geq2$, let $\hat{P}=(1-\sum_{n=2}^{\infty}p_n)I_{\omega=1}+\sum_{n\geq2}p_nI_{\omega=n}$. The optimal estimator
\begin{equation}
\begin{aligned}
\hat{\eta}=E_{\hat{P}}[\xi|\mathcal{C}]=E_{\hat{P}}[\xi]=(1-\sum_{n=2}^{\infty}p_n)+\sum_{n=2}^{\infty}\frac{2^n}{n^4}p_n=1+\sum_{n=2}^{\infty}(\frac{2^n}{n^4}-1)p_n,
\end{aligned}
\end{equation}
and
\begin{equation}
\begin{aligned}
E_{P}[\xi-\hat{\eta}]^2=[\frac{\lambda}{2}+\frac{2(1-\lambda)}{3}][\sum_{n=2}^{\infty}(\frac{2^n}{n^4}-1)p_n]^2+\sum_{n=2}^{\infty}[\frac{\lambda}{2^n}+\frac{2(1-\lambda)}{3^n}][(\frac{2^n}{n^4}-1)-\sum_{m=2}^{\infty}(\frac{2^m}{m^4}-1)p_m]^2.
\end{aligned}
\end{equation}
By the optimal conditions $\frac{\partial E_P[\xi-\hat{\eta}]^2}{\partial p_i}=0$ and $\frac{\partial E_{P}[\xi-\hat{\eta}]^2}{\partial \lambda}=0$, we deduce that
\begin{equation}\label{p value}
\left\{
\begin{aligned}
p_i&=\frac{\lambda}{2^i}+\frac{2(1-\lambda)}{3^i},\quad i\geq2,\\
\lambda&=\frac{F(n)-\sum_{m=2}^{\infty}\frac{2}{3^m}(\frac{2^m}{m^4}-1)}{\sum_{m=2}^{\infty}(\frac{2^m}{m^4}-1)(\frac{1}{2^m}-\frac{2}{3^m})}
\end{aligned}
\right.
\end{equation}
where $F(n)=\frac{\sum_{n=2}^{\infty}(\frac{1}{2^n}-\frac{2}{3^n})(\frac{2^n}{n^4}-1)^2}{2\sum_{n=2}^{\infty}(\frac{1}{2^n}-\frac{2}{3^n})(\frac{2^n}{n^4}-1)}.$\\
\rightline{$\square$}
\begin{example}
Given a complete filtered probability space $(\Omega, \mathcal{F}, \{\mathcal{F}_t\}_{0\leq t\leq T},P_0)$, $W(\cdot)$ is a standard one dimensional Brownian motion defined on this space where $\mathcal{F}_t=\sigma\{W(s),0\leq s\leq t\}$ and $\mathcal{F}=\mathcal{F}_T$. The space $L_{\mathcal{F}}^{4+2\epsilon}(0,T;\mathbb{R})$ denotes all the $\mathcal{F}_t$-progressively measurable processes $h_t$ such that $E_{P_0}\int^T_0|h_t|^{4+2\epsilon}dt<\infty$ for given constant $\epsilon\in (0,1)$. Let us introduce g-expectation defined by the following backward stochastic differential equation:
\begin{equation}\label{bsde}
y_t=\xi+\int^T_0|z_s|ds+\int^T_0z_sdW(s)
\end{equation}
where $\xi$ is a $\mathcal{F}_T$-measurable $(4+2\epsilon)$ integrable random variable. Here $g(y,z)=|z|$. According to the results in \cite{PP}, there exists a unique adapted pair $\{y_t,z_t\}_{t\in{0,T}}$ which solves \eqref{bsde}. We call the solution $\{y_t\}_{0\leq t\leq T}$ the conditional g-expectation with respect to $\mathcal{F}_t$ and denote it by $\mathcal{E}_{|z|}(\xi|\mathcal{F}_t)$.

Consider the following linear case:
\begin{equation}\label{lbsde}
\tilde{y}_t=\xi+\int^T_t\mu_sz_sds+\int^T_tz_sdW(s)
\end{equation}
where $|\mu_s|\leq 1,$ $P_0-a.s.$. By Girsanov transform, there exists a probability $P^{\mu}$ such that $\{y_t\}_{0\leq t \leq T}$ of \eqref{lbsde} is a martingale under $P^{\mu}$. Let $\mathcal{P}:=\{P^{\mu}\big||\mu_s|\leq 1, \ P_0-a.s.\}$. By Theorem 2.1 in \cite{CJ},
$$\mathcal{E}_{|z|}(\xi)=\sup_{P^{\mu}\in \mathcal{P}}E_{P^{\mu}}[\xi], \quad \forall \xi \in L_{\mathcal{F}_T}^{4+2\epsilon}(P_0) $$
and
$$\mathcal{E}_{|z|}(\xi|\mathcal{F}_t)=\esssup_{P^{\mu}\in \mathcal{P}}E_{P^{\mu}}[\xi|\mathcal{F}_t], \quad \forall \xi \in L_{\mathcal{F}_T}^{4+2\epsilon}(P_0). $$

It is easy to see that $\mathcal{E}_{|z|}(\cdot)$ is a sublinear operator. Denote the corresponding minimum mean square estimator by $\rho_{|z|}(\xi|\mathcal{F}_t)$. We claim that the minimum mean square estimator $\rho_{|z|}(\xi|\mathcal{F}_t)$ does not coincide with $\mathcal{E}_{|z|}(\xi|\mathcal{F}_t)$. Otherwise, If not, i.e.$\rho_{|z|}(\xi|\mathcal{F}_t)=\mathcal{E}_{|z|}(\xi|\mathcal{F}_t)$, as the result of Proposition \ref{proposition}, we have
$$\esssup_{P^{\mu}\in \mathcal{P}}E_{P^{\mu}}[\xi|\mathcal{F}_t]=\rho_{|z|}(\xi|\mathcal{F}_t)=-\rho_{|z|}(-\xi|\mathcal{F}_t)=\essinf_{P^{\mu}\in \mathcal{P}}E_{P^{\mu}}[\xi|\mathcal{F}_t].$$
Since $\mathcal{P}$ contains more than one probability measures, the above equation can not be true for all the $(4+2\epsilon)$ integrable $\xi\in\mathcal{F}_T$. Thus, our claim holds.\\
\rightline{$\square$}
\end{example}

\section{Acknowledgment}

The authors would like to thank editors and an anonymous referee for helpful comments and suggestions, which lead to a much better
version of this paper.

\begin{appendices}
\section{}
In this section, we give the following lemma which is used to prove Theorem \ref{unique}.
\begin{lemma}\label{calculate}
Let $\hat{\eta}_1=E_{\hat{P}_1}[\xi|\mathcal{C}]$, $\hat{\eta}_2=E_{\hat{P}_2}[\xi|\mathcal{C}]$, $P^{\lambda}=\lambda\hat{P}_1+(1-\lambda)\hat{P}_2$, $\lambda_{\hat{P}_1}=\lambda E_{P^{\lambda}}\big[\frac{d\hat{P}_1}{dP^{\lambda}}\big]$, $\lambda_{\hat{P}_2}=(1-\lambda) E_{P^{\lambda}}\big[\frac{d\hat{P}_2}{dP^{\lambda}}\big]$. Then we have
$$E_{P^{\lambda}}[(\xi-\lambda_{\hat{P}_1}\hat{\eta}_1-\lambda_{\hat{P}_2}\hat{\eta}_2)^2]\geq \alpha.$$
\end{lemma}
$\textbf{Proof}$.
\begin{equation}\label{A1}
\begin{aligned}
&E_{P^{\lambda}}[(\xi-\lambda_{\hat{P}_1}\hat{\eta}_1-\lambda_{\hat{P}_2}\hat{\eta}_2)^2]\\
=&E_{P^{\lambda}}\big[\big(\lambda_{\hat{P}_1}(\xi-\hat{\eta}_1)+\lambda_{\hat{P}_2}(\xi-\hat{\eta}_2)\big)^2\big]\\
=&E_{P^{\lambda}}\big[\lambda^2_{\hat{P}_1}(\xi-\hat{\eta}_1)^2+\lambda^2_{\hat{P}_2}(\xi-\hat{\eta}_2)^2+2\lambda_{\hat{P}_1}\lambda_{\hat{P}_2}(\xi-\hat{\eta}_1)(\xi-\hat{\eta}_1)\big]\\
=&E_{P^{\lambda}}\big[\lambda_{\hat{P}_1}(\xi-\hat{\eta}_1)^2+\lambda_{\hat{P}_2}(\xi-\hat{\eta}_2)^2-\lambda_{\hat{P}_1}\lambda_{\hat{P}_2}(\hat{\eta}_1-\hat{\eta}_2)^2\big]\\
=&\lambda E_{\hat{P}_1}[\lambda_{\hat{P}_1}(\xi-\hat{\eta}_1)^2]+(1-\lambda)E_{\hat{P}_2}[\lambda_{\hat{P}_1}(\xi-\hat{\eta}_1)^2]+\lambda E_{\hat{P}_1}[\lambda_{\hat{P}_2}(\xi-\hat{\eta}_2)^2]
+(1-\lambda)E_{\hat{P}_2}[\lambda_{\hat{P}_2}(\xi-\hat{\eta}_2)^2]\\&-\lambda E_{\hat{P}_1}[\lambda_{\hat{P}_1}\lambda_{\hat{P}_2}(\hat{\eta}_1-\hat{\eta}_2)^2]-(1-\lambda)E_{\hat{P}_2}[\lambda_{\hat{P}_1}\lambda_{\hat{P}_2}(\hat{\eta}_1-\hat{\eta}_2)^2]\\
=&\lambda E_{\hat{P}_1}[(\xi-\hat{\eta}_1)^2]-\lambda E_{\hat{P}_1}[\lambda_{\hat{P}_2}(\xi-\hat{\eta}_1)^2]+(1-\lambda)E_{\hat{P}_2}[(\xi-\hat{\eta}_1)^2]-(1-\lambda)E_{\hat{P}_2}[\lambda_{\hat{P}_2}(\xi-\hat{\eta}_1)^2]\\
&+\lambda E_{\hat{P}_1}[(\xi-\hat{\eta}_2)^2]-\lambda E_{\hat{P}_1}[\lambda_{\hat{P}_1}(\xi-\hat{\eta}_2)^2]+(1-\lambda) E_{\hat{P}_2}[(\xi-\hat{\eta}_2)^2]-(1-\lambda) E_{\hat{P}_2}[\lambda_{\hat{P}_1}(\xi-\hat{\eta}_2)^2]\\
&-\lambda E_{\hat{P}_1}[\lambda_{\hat{P}_1}\lambda_{\hat{P}_2}(\hat{\eta}_1-\hat{\eta}_2)^2]-(1-\lambda)E_{\hat{P}_2}[\lambda_{\hat{P}_1}\lambda_{\hat{P}_2}(\hat{\eta}_1-\hat{\eta}_2)^2].
\end{aligned}
\end{equation}

Because
$$(1-\lambda)E_{\hat{P}_2}[(\xi-\hat{\eta}_1)^2]=(1-\lambda)E_{\hat{P}_2}[(\lambda_{\hat{P}_1}+\lambda_{\hat{P}_2})(\xi-\hat{\eta}_1)^2]$$
and
$$\lambda E_{\hat{P}_1}[(\xi-\hat{\eta}_2)^2]=\lambda  E_{\hat{P}_1}[(\lambda_{\hat{P}_1}+\lambda_{\hat{P}_2})(\xi-\hat{\eta}_2)^2],$$
then equation (6.1) becomes

\begin{equation*}
\begin{aligned}
(6.1)&=\lambda E_{\hat{P}_1}[\lambda_{\hat{P}_2}(\xi-\hat{\eta}_2)^2-\lambda_{\hat{P}_2}(\xi-\hat{\eta}_1)^2]+(1-\lambda)E_{\hat{P}_2}[\lambda_{\hat{P}_1}(\xi-\hat{\eta}_1)^2-\lambda_{\hat{P}_1}(\xi-\hat{\eta}_2)^2]\\
&-\lambda E_{\hat{P}_1}[\lambda_{\hat{P}_1}\lambda_{\hat{P}_2}(\hat{\eta}_1-\hat{\eta}_2)^2]-(1-\lambda)E_{\hat{P}_2}[\lambda_{\hat{P}_1}\lambda_{\hat{P}_2}(\hat{\eta}_1-\hat{\eta}_2)^2]+\lambda E_{\hat{P}_1}[(\xi-\hat{\eta}_1)^2]\\
&+(1-\lambda)E_{\hat{P}_2}[(\xi-\hat{\eta}_2)^2].
\end{aligned}
\end{equation*}

Firstly, we calculate the items in the expectation operator $\lambda E_{\hat{P}_1}[\cdot]$
\begin{equation*}
\begin{aligned}
&\lambda_{\hat{P}_2}(\xi^2-\hat{\eta}_2^2-2\xi\hat{\eta}_2)-\lambda_{\hat{P}_2}(\xi^2-\hat{\eta}_1^2-2\xi\hat{\eta}_1)-\lambda_{\hat{P}_1}\lambda_{\hat{P}_2}(\hat{\eta}_1-\hat{\eta}_2)^2\\
=&\lambda_{\hat{P}_2}[2\hat{\eta}_1^2(\hat{\eta}_2-\hat{\eta}_1)+2\xi(\hat{\eta}_1-\hat{\eta}_2)]+\lambda_{\hat{P}_2}^2(\hat{\eta}_1-\hat{\eta}_2)^2\\
=&\lambda_{\hat{P}_2}[2(\xi-\hat{\eta}_1)(\hat{\eta}_1-\hat{\eta}_2)]+\lambda_{\hat{P}_2}^2(\hat{\eta}_1-\hat{\eta}_2)^2.
\end{aligned}
\end{equation*}

Because $\lambda_{\hat{P}_2}(\hat{\eta}_1-\hat{\eta}_2)$ is $\mathcal{C}-$ measurable and $(\xi-\hat{\eta}_1)$ is orthogonal with $\sigma$- algebra $\mathcal{C}$ under probability measure $\hat{P}_1$, it results that
$$\lambda E_{\hat{P}_1}[\lambda_{\hat{P}_2}2(\xi-\hat{\eta}_1)(\hat{\eta}_1-\hat{\eta}_2)]=\lambda E_{\hat{P}_1}[\lambda_{\hat{P}_2}(\hat{\eta}_1-\hat{\eta}_2)]E_{\hat{P}_1}[2(\xi-\hat{\eta}_1)]=0.$$

Similarly, we can also calculate the items in the expectation operator $(1-\lambda) E_{\hat{P}_2}[\cdot]$. Then equation \eqref{A1}  becomes

\begin{equation*}
\begin{aligned}
&E_{P^{\lambda}}[(\xi-\lambda_{\hat{P}_1}\hat{\eta}_1-\lambda_{\hat{P}_2}\hat{\eta}_2)^2]\\
=&\lambda E_{\hat{P}_1}\big[(\xi-\hat{\eta}_1)^2\big]+(1-\lambda) E_{\hat{P}_2}\big[(\xi-\hat{\eta}_2)^2\big]\\
&+\lambda E_{\hat{P}_1}\big[\lambda^2_{\hat{P}_2}(\hat{\eta}_1-\hat{\eta}_2)^2\big]+(1-\lambda)E_{\hat{P}_2}\big[\lambda^2_{\hat{\eta}_1}(\hat{\eta}_1-\hat{\eta}_2)^2\big].
\end{aligned}
\end{equation*}
\rightline{$\square$}

\begin{definition}
For a given probability space $(\Omega,\mathcal{F},P_0)$, $\{\mathcal{F}_n\}_{n\geq1}$ is the filtration satisfying $\mathcal{F}:=\bigvee_{n=1}\mathcal{F}_n$. We say that the set $\mathcal{P}$ is stable if for elements $Q^{0},\ Q\in \mathcal{P}^{e}$ with associated martingales $Z_n^{0},\ Z_n$ and for each stopping time $\tau$, the martingale L defined as $L_n=Z_n^0$ for $n\leq \tau$ and $L_n=Z_{\tau}^0\frac{Z_n}{Z_{\tau}}$ for $n\geq \tau$ defines an element of $\mathcal{P}$, where $\mathcal{P}^e$ denotes the elements in $\mathcal{P}$ which is equivalent to $P_0$ and $Z_n^{Q}:=E_{P_0}[\frac{dQ}{dP_0}|\mathcal{F}_n]$.
\end{definition}
\end{appendices}
\renewcommand{\refname}

\end{document}